\documentclass{article}[12pt]
\usepackage[letterpaper,margin=1in]{geometry}
\usepackage{amsmath,amssymb,graphicx,subcaption}
\usepackage{color}

\newcommand \R{\mathbb{R}}

\usepackage{url}

\title{Degenerate Umbilic Points of Analytic Surfaces}
\author{John Guckenheimer\\
Mathematics Department\\
Cornell University}

\begin{document}
\maketitle
\begin{abstract}

Umbilics are points of a surface embedded in three space where normal curvatures are independent of direction. The (in)famous Carath\'{e}odory Conjecture states that a compact simply connected embedded surface has at least two umbilic points. A counterexample to this conjecture would be a surface whose principal foliations have index two at a single umbilic. All (purported) proofs of the Carath\'{e}odory Conjecture are based on analyses of the index of an umbilic, concluding that it is at most one. This paper investigates principal foliations near degenerate umbilics, giving a much simpler geometric argument than previous work that the index of an umbilic on an analytic surface cannot be an integer larger than one. The results also establish lower bounds for the index of an umbilic based on its degeneracy.

\end{abstract}

\section*{Principal Curvatures and Foliations}

Let $S \subset \R^3$ be a smooth oriented surface. If $p \in S$, $N(p)$ is the unit normal to $S$ at $p$, $\gamma(s)$ is a planar curve on $S$ parameterized by arc length with $\gamma(0) = p$ and $\frac{d\gamma}{ds}(0) = v$, then the normal curvature at $p$ in the direction $v$ is $\kappa(p,v) = \frac{d^2(\gamma)}{ds^2}(0) \cdot N(p)$. \emph{Principal vectors} realize the maximum and minimum normal curvatures as a function of direction $v$ at each point of $S$. \emph{Rodrigues formula}~\cite{MR2007065} states that the tangent vector $\frac{d\gamma}{dt}(p)$ to $\gamma$ is a principal vector $v$ if and only if $ \frac{dN \circ \gamma}{dt}(p) = \kappa \frac{d\gamma}{dt}(p)$ for some $\kappa$, producing a system of equations whose solutions $v$,$\kappa$ are principal unit vectors:
\begin{equation}
\begin{split}
 N \cdot v & = 0 \\
 v \cdot v & = 1 \\
 \det(N,v,dN \cdot v) & = 0. 
\end{split} 
\label{rod}
\end{equation}
The equations \eqref{rod} express that $v$ is a tangent vector to $S$, $v$ has unit length and that $dN \cdot v$ lies in the plane spanned by $N$ and $v$. The first of these equations is linear in $v$ while the second two are quadratic and even in $v$. Points where the normal curvatures are independent of direction are \emph{umbilics}(or umbilic points). In the complement of the umbilics, the Rodrigues formula defines a \emph{cross field} consisting of two orthogonal line fields aligned with the principal vectors. Integration of these line fields yields the \emph{principal foliations} of the surface with singularities at the umbilics. The \emph{index} of an isolated umbilic measures the change in the angle of a continuous curve of principal vectors along a positively oriented loop surrounding that umbilic and no others. Umbilics of generic surfaces were classified by Darboux~\cite{Darboux}:  they are now called \emph{lemons}, \emph{monstars} and \emph{stars} with indices $ \frac{1}{2},\frac{1}{2}$ and $-\frac{1}{2}$~\cite{Berry_1977}. Principal foliations of generic surfaces are not orientable since continuation of a principal vector around a loop containing a single generic umbilic reverses its orientation. On a simply connected compact surface with isolated umbilics, the index theorem concludes the the sum of their indices is two. 

Beyond bifurcations in generic one parameter families of principal foliations~\cite{MR2105778}, existing literature on degenerate umbilics is dominated by the Carath\'{e}odory Conjecture. This century old conjecture states that a compact, simply connected surface in $\R^3$ has at least two umbilics. The status of the Carath\'{e}odory Conjecture is controversial, with papers claiming to prove the conjecture followed by work describing deficiencies of the proofs. The goal of most work on the Carath\'{e}odory conjecture has been to establish that the index of an isolated umbilic point is at most one. Since the sum of the indices of umbilics of a compact simply connected surface is two, this local result implies the Carath\'{e}odory conjecture. This paper presents a new analysis of how the indices of umbilics change as they merge with one another in analytic surfaces.

\section*{Monge Coordinates at Umbilics}

Denote the set of umbilics on the surface $S$ by $\Upsilon$. If $ \chi \in \Upsilon$ is an umbilic, choose orthonormal coordinates $(x,y,z)$ so that $\chi$ is the origin and the $z$ axis is normal to $S$. The implicit function theorem implies that there is a neighborhood $U$ of $\chi$ in which $S \cap U$ is the graph of a function $z = h(x,y)$ with $dh_{(0,0)} = (0,0)$. These are \emph{Monge coordinates}.
Rodrigues formula gives an equation for the principal directions of S in terms of the derivatives of $h$. Denote the first and second partial derivatives of $h$ by $h_x, h_y, h_{xx}, h_{xy}, h_{yy}$. The unit vector field $N = \alpha (-h_x,-h_y,1)^\mathsf{T}, \; \alpha^{-2} = 1 + h_x^2+h_y^2$ is normal to $S$. Differentiating, 
\begin{equation}
dN = d\alpha \begin{pmatrix} -h_x \\ -h_y \\1 \end{pmatrix} + \alpha \begin{pmatrix} -h_{xx} & -h_{xy} & 0 \\ -h_{xy} & -h_{yy} & 0 \\ 0 & 0 & 0   \end{pmatrix}.
\end{equation}
A non-zero vector $(u,v,w)^\mathsf{T}$ is tangent to $S$ if and only if $w = u h_x + v h_y$.
Rodrigues Formula implies that it is a principal vector when
\begin{equation}
R(x,y,u,v,w) = h_{xy} h_y v w - h_x h_{yy} v w + h_{xy} v^2 + h_{xx}h_y u w - h_x h_{xy} u w - h_{yy} u v + h_{xx} u v - h_{xy} u^2 
\label{rodh}
\end{equation}
vanishes with $w = u h_x + v h_y$.
Equation~\eqref{rodh} is satisfied identically when $S$ is a sphere. Spheres and planes are the only surfaces all of whose points are umbilics~\cite{MR3837152}. Analytic Monge surfaces with an isolated umbilic at the origin
are graphs of functions 
\begin{equation}
h(x,y) = 1-\sqrt{1-x^2-y^2} + p(x,y) + q(x,y) 
\label{spert}
\end{equation}
where $p$ is a non-zero homogeneous polynomial of degree $d$ and $q$ is analytic function that is $o(d)$ at the origin. 
We observe that $R(x,y,u,v,u h_x + v h_y)$ is a homogeneous quadratic polynomial of $(u,v)$ for each $(x,y)$. Writing $R(x,y,u,v,u h_x + v h_y) = c_{uu} u^2 + c_{uv} u v + c_{vv} v^2$, and $p_a(x,y) = \sum_{k=0}^d a_k x^{d-k}y^k, \; a = (a_0, \cdots, a_d)$, we find with the aid of symbolic computations in Matlab~\cite{MATLAB2024b} that
\begin{equation}
\begin{split}
c_{uu} = & -\sum_{k = 1}^{d-1} k(d-k) a_{k} x^{d-1-k}y^{k-1}  + O(d-1)\\
c_{uv} = & \sum_{k = 1}^{d-1} ((d-k)(d-k-1) a_{k-1} - k(k+1) a_{k+1})x^{d-1-k}y^{k-1} + O(d-1)\\
c_{vv} = & \sum_{k = 1}^{d-1} k(d-k) a_{k} x^{d-1-k}y^{k-1} + O(d-1)
\label{rodc}
\end{split}
\end{equation}
Note that $c_{uu} + c_{vv} = O(d-1)$. 
On circles $C_r, \, r^2 = x^2+y^2$, centered at the origin in the $(x,y)$ plane, the equations \eqref{rodh} scaled by $r^{2-d}$ have a limit as $x^2+y^2 \to 0$ that depends only on the degree $(d-2)$ terms $c_{uu},c_{uv},c_{vv}$ in $(x,y)$ displayed in equations \eqref{rodc} .  We denote these degree $(d-2)$ terms by $\bar{C}_{uu},\bar{C}_{uv},\bar{C}_{vv}$. When they do not vanish simultaneously at a point on $C_r$, they determine the index of the umbilic at the origin. If $d$ is odd, then the index is an odd multiple of $1/2$ since following solutions of $R=0$ continuously around $C_r$ reverses their orientation. Thus, umbilic points of index 2 do not occur for odd $d$. 
Consequently, $S_a^d$ that are graphs of $h(x,y) = 1-\sqrt{1-x^2-y^2} + p_a(x,y)$ with $d>2$ even are the ones most relevant to the Carath\'{e}odory conjecture. 

\section*{Umbilic Indices}

This section investigates the index $I_a(0)$ of the origin in the family of surfaces $S_a^d$ with $d>3$ an even integer. The group $SO (2)$ of rotations of the $(x,y)$ plane acts on the family $S^d_a$, preserving the index of the origin. If $\chi_a$ is a non-zero umbilic, there is a rotation that carries $\chi_a$ to the $x-$axis. The value $c_{uu}u^2 + c_{uv}uv + c_{vv}v^2 = R(x,0,u,v,u h_x + v h_y)$ of Rodrigues formula \eqref{rodh} at $(x,0)$ depends only on the coefficients $(a_0,a_1,a_2)$ because all the terms in the derivatives $h_x,h_y,h_{xx},h_{xy},h_{yy}$ with coefficients $a_k, k>2$ are divisible by $y$. Furthermore, the terms $c_{uu}$ and $c_{vv}$ are divisible by $a_1x^{d-2}$ while 
$$c_{uv} = (d(d-1)a_0 - 2a_2)x^{d-2} - \frac{2(da_0+a_2)x^d}{1-x^2} + O(x^{2d-2}).$$
These facts imply that $S_a$ has umbilics on the $x-$axis near the submanifold of the parameter space defined by $a_1 = 0, \;  d(d-1)a_0 - 2a_2 - 2(da_0+a_2)x^2 = 0$. In a family with umbilics on the $x$-axis, these umbilics merge with the origin when $d(d-1)a_0 = 2a_2$. 

There are some parameter values $a$ for which the index of the origin is easily determined from symmetry considerations. First, the graph of $h(x,y) =  1-\sqrt{1-x^2-y^2} + (x^2+y^2)^{d/2}$ is a surface of revolution, so its principal foliations consist of circles centered at the origin and rays emanating from the origin. The index is $+1$. Second, the polynomial $((x+iy)^d + (x-iy)^d)/2$ is invariant under rotations of the $(x,y)$ plane by angles that are integer multiples of  $2\pi/d$. Thus the graphs of functions $h(x,y) =  1-\sqrt{1-x^2-y^2} + \alpha_1(x^2+y^2)^{d/2} + \alpha_2((x+iy)^d + (x-iy)^d))/2$ are also invariant under rotation around the $z$ axis by $2\pi/d$. The following result shows that $I_0 \ge 1-d/2$.

\vspace{6pt}
\noindent
\textbf{Theorem:} In the family $S^d_a$, $1-d/2 \le I_a(0)$ whenever $(\bar{C}_{uu},\bar{C}_{uv})$ does not vanish. 

\vspace{6pt}
\noindent
\textbf{Proof:} Fix $a$ and examine $S_a$, the graph of the function
\begin{equation*}
 h_a(x,y) = 1-\sqrt{1-x^2-y^2} + \sum_{j=0}^d a_j x^{(d-j)}y^j.
\end{equation*}
The Rodrigues formula $R(x,y,u,v,uh_x+vh_y)$ is a homogeneous polynomial in $(u,v)$ and an analytic function in $(x,y)$ whose lowest degree terms have degree $(d-2)$. Thus the truncation $\bar{R}$ of $R$ to terms of degree $d-2$ can be viewed as a polynomial equation in projective coordinates on $\mathbb{P} \times \mathbb{P}$ of degree $2$ in the $(u,v)$ coordinate and degree $d-2$ in the $(x,y)$ coordinate. When the solutions to this equation  form a smooth curve, $I_a(0)$ is defined and determined by the curve. The curve has winding number $2$ along the $(u,v)$ coordinate and winding number in the interval $[2-d,d-2]$ along the $(x,y)$ coordinate. Geometrically, the solutions represent the winding of the principal foliations of $S$ along a infinitesimal curve that encircles the origin twice. This implies that $1-d/2 \le I_a(0) \le d/2-1$. 

\vspace{6pt}
The main result of this paper is the following:

\vspace{6pt}
\noindent
\textbf{Theorem:}
 Let $S_{a_\epsilon} \subset S_a^d, \epsilon \in [0,1]$ be a one parameter family of surfaces with $I_0(S_{a_\epsilon}) = 1$ for $\epsilon < \epsilon_\ast$ but $I_0(S_{a_\epsilon}) \ne 1$ for $\epsilon > \epsilon_\ast$. Then $I_0(S_{a_\epsilon}) < 1$ for $\epsilon > \epsilon_\ast$. Generic umbilics $\chi_\epsilon$ emerging from the origin at $\epsilon = \epsilon_\ast$ are stars.

\vspace{6pt}
\noindent
\textbf{Proof:}
Berry and Hannay~\cite{Berry_1977} write the Monge form of a surface at an umbilic as
$$f(x,y) = \frac{1}{2}k(x^2+y^2)+\frac{1}{6}(\alpha x^3 + 3\beta x^2y + 3\gamma xy^2 + \delta y^3) + O(4).$$
They further show that the index of generic umbilics is determined by the quantity 
$$J(\alpha,\beta,\gamma,\delta)  = \alpha \gamma - \gamma^2 + \beta \delta - \beta^2.$$ 
If $J$ is positive, the umbilic index is $\frac{1}{2}$, while if $J$ is negative, the umbilic index is $-\frac{1}{2}$ and the umbilic is a star. Thus to prove the theorem, it suffices to establish that $J$ is negative at a non-zero umbilic $\chi_\epsilon$ with $\epsilon$ slightly larger than $\epsilon_\ast$. 

Without loss of generality, consider the surface $S_{a}$ which is the graph of $ h_a(x,y) = 1-\sqrt{1-x^2-y^2} + \sum_{j=0}^d a_j x^{(d-j)}y^j$ with $a = a_\epsilon$ and umbilic $\chi = (x_u,0,H_a(x_u,0,0))$. Set
$$H_a(x,y,z) = z - (1-\sqrt{1-x^2-y^2} + \sum_{j=0}^d a_j x^{(d-j)}y^j)$$ 
to be a function of $(x,y,z)$ whose zero level set is $S_a$.
Compute the Monge form $\hat{h}$ of the surface $S_a$ at $\chi$ in three steps:
\begin{itemize}
\item
 Set $H^1(x,y,z) = H^1(x-x_u,y,z- h_a(x_u,0,0))$ to translate $\chi$ to the origin.
 \item
 Set $\tilde{H} = H^1 \circ T$ where T is the rotation of the $(x,z)$ plane so that the normal $d\tilde{H}(0,0,0) = d(H^1 \circ T)(0,0,0) = (0,0,1)$ is in the vertical direction. Define new coordinates by $(\hat{x},\hat{y},\hat{z}) = T(x-x_u,y,z-h_a(x_u,0,0))$.
 \item
Replace $\tilde{H}$ by a function $\hat{H}$ that has the same zero level set, but is in Monge form: $\hat{H}(\hat{x},\hat{y},\hat{z}) = \hat{z} - \hat{h}(\hat{x},\hat{y})$. 
The function $\hat{H}$ will have the form $g(\hat{x},\hat{y},\hat{z})\tilde{H}(\hat{x},\hat{y},\hat{z})$ with $g$ determined by solving for Taylor series coefficients that put $g\tilde{H}$ into the desired form $\hat{z} - \hat{h}(\hat{x},\hat{y})$. 
The following terms of degrees 1, 2 and 3 may have non-zero coefficients in the Taylor series of $\tilde{H}$: $\hat{z}, \hat{x}^2, \hat{y}^2, \hat{x}\hat{z}, \hat{z}^2, \hat{x}^3,  \hat{x}^2\hat{z}, \hat{x}\hat{y}^2, \hat{x}\hat{z}^2, \hat{y}^2\hat{z}, \hat{z}^3$. 
Setting $g(\hat{x},\hat{y},\hat{z}) = b_x\hat{x} + b_z\hat{z} + b_{xx}\hat{x}^2 + b_{xz}\hat{x}\hat{z} + b_{zz}\hat{z}^2 + b_{yy}\hat{y}^2$ with 
$$b_x = -\hat{H}_{xz}/\hat{H}_z, b_z = -\hat{H}_{zz}/\hat{H}_z, b_{xx} = -\hat{H}_{xxz}/\hat{H}_z +\hat{H}_{xz}^2/\hat{H}_z^2 + \hat{H}_{xx}\hat{H}_{zz}/\hat{H}_z^2, $$
$$b_{xz} =  -\hat{H}_{xzz}/\hat{H}_z + 2\hat{H}_{zz}\hat{H}_{xz}/\hat{H}_z^2, b_{zz} = -\hat{H}_{zzz}/\hat{H}_z +\hat{H}_{zz}^2/\hat{H}_z^2, b_{yy} = -\hat{H}_{yyz}/\hat{H}_z + \hat{H}_{yy}\hat{H}_{zz}/\hat{H}_z^2$$
 eliminates the terms $\hat{x}\hat{z},\hat{z}^2,\hat{x}^2\hat{z},\hat{x}\hat{z}^2,\hat{z}^3$  from $\hat{H}$ and 
 replaces the coefficients of the terms $\hat{H}_{xxx},\hat{H}_{xyy}$ by $\hat{H}_{xxx} - \hat{H}_{xx} \hat{H}_{zz} / \hat{H}_z, \hat{H}_{xyy} - \hat{H}_{xz} \hat{H}_{yy}/\hat{H}_z$. 
 Finally, divide $\hat{H}$ by $\hat{H}_z$ to obtain the desired Monge form $H$.
 \end{itemize}
 The result of applying this computation to $H_a(x,y,z) = z - (1-\sqrt{1-x^2-y^2} +  a_0 x^d + a_2 x^{d-2} y^2)$ yields to leading order in $x_u$ the Monge form 
 $$H(\hat{x},\hat{y},\hat{z}) = \frac{1}{2}(\hat{x}^2+\hat{y}^2) - a_0 x_u^{d-3} d(d-1)(d-2)(\frac{1}{6}\hat{x}^3 + \frac{1}{2}\hat{x}\hat{y}^2)$$
 The value of $J$  to leading order in $x_u$ is $ -d^2(d^2-1)(d-2)/2$ which is clearly negative for $d>2$. Consequently, any generic umbilic close to the origin is a star with index $-\frac{1}{2}$. Finally,  consider the merger of umbilics $\chi_\epsilon$ with the origin at $\epsilon_\ast$ as $\epsilon$ varies. For $\epsilon_0$ slightly smaller than $\epsilon_\ast$, choose a circle $C_r$ centered at the origin which contains no other umbilics. Then its winding number is $1$, the index of the origin. As $\epsilon < \epsilon_\ast$ increases, one or more umbilics cross $C_r$. These umbilics are stars, so the winding number of $C_r$ becomes smaller than $1$. As these umbilic merge with the origin at $\epsilon_\ast$, they decrease its index. Thus, the index of the origin is smaller than $1$ when $\epsilon > \epsilon_\ast$. The theorem is proved.

\section*{Visualization of Principal Foliations in $S^8_a$}

This section computes the index $I_0(a)$ of the origin and displays the principal foliations of a few surfaces in $S^8_a$. The Rodrigues formula $R(x,y,u,v,uh_x+vh_y) = R(x,y,u,v)$ for the surface 
\begin{equation}
 h_a(x,y) = 1-\sqrt{1-x^2-y^2} + \sum_{j=0}^8 a_j x^{(8-j)}y^j.
\end{equation} 
in $S^8_a$ is
\begin{equation}
 \begin{split}
 R(x,y,u,v) =  
 (- 7 a_7 y^6 - 12 a_6 x y^5 - 15 a_5 x^2 y^4 - 16 a_4 x^3 y^3 - 15 a_3 x^4 y^2 - 12 a_2 x^5 y - 7 a_1 x^6) & u^2 \\
 +((2 a_6 - 56 a_8) y^6 + (6 a_5 - 42 a_7) x y^5 + (12 a_4 - 30 a_6) x^2 y^4 + (20 a_3 - 20 a_5) x^3 y^3 \\
 + (30 a_2 - 12 a_4) x^4 y^2 + (42 a_1 - 6 a_3) x^5 y + (56 a_0 - 2 a_2) x^6) & uv \\
 +( 7 a_7 y^6 + 12 a_6 x y^5  + 15 a_5 x^2 y^4  + 16 a_4 x^3 y^3  + 15 a_3 x^4 y^2  + 12 a_2 x^5 y  + 7 a_1 x^6) & v^2 + R_7
 \end{split}
 \label{htay}
\end{equation}
where the minimal $(x,y)$ degree of terms in $R_7$ is 7. Writing $R(x,y,u,v,uh_x+vh_y) = c_{uu}(x,y,a)u^2 + c_{uv}(x,y,a)uv + c_{vv}(x,y,a)v^2 $, $(x,y,a)$ is an umbilic point when $c_{uu} = c_{uv} = c_{vv} = 0$. Note that the sum of the coefficients of $u^2$ and $v^2$ is $O(7)$ in $(x,y)$. 

\begin{figure}[ht]
\centering
\includegraphics[height=1.75in]{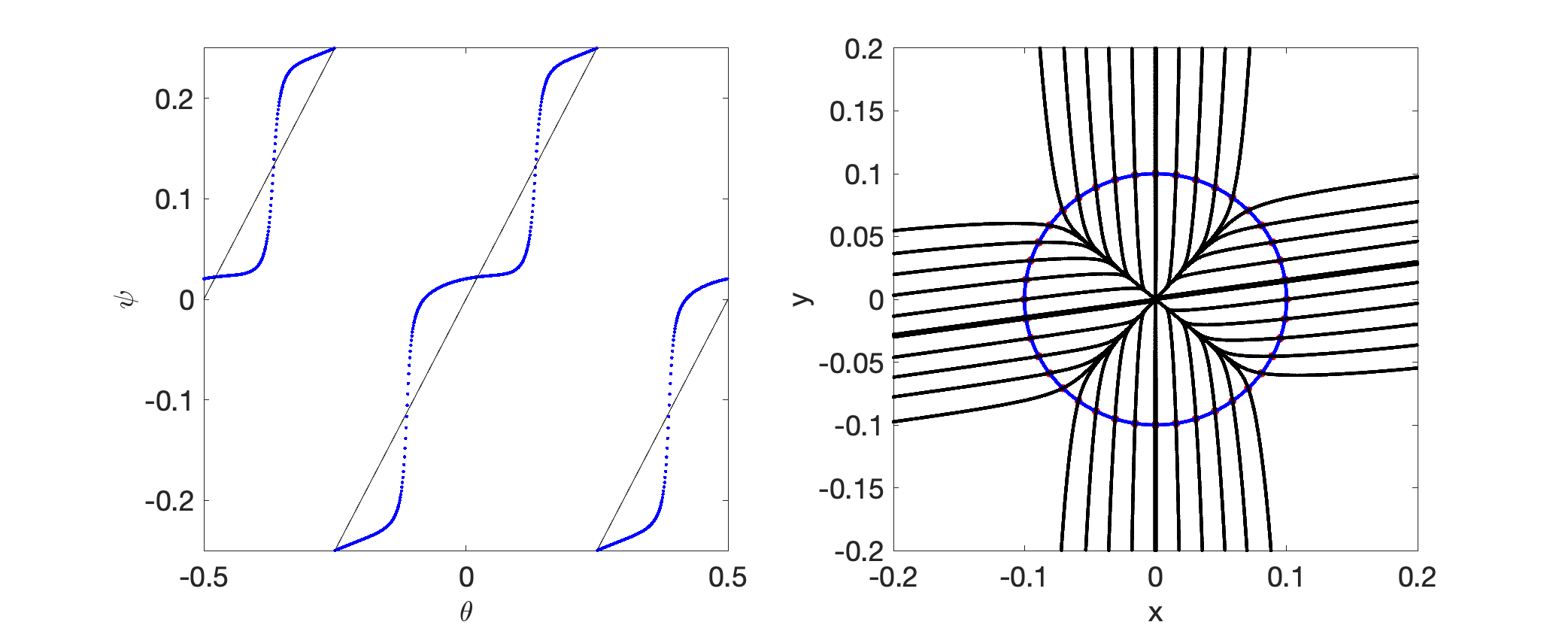}
  \caption{The index $I_0(a)$ for parameter $a = (1,1,1,0,0,0,1,0,1)$. The left panel shows 1000 solutions of $R(0.1 \cos(2\pi\theta),0.1 \sin(2\pi\theta),\cos(2\pi\psi),\sin(2\pi\psi)) = 0$. The full set of solutions is a curve that winds once around both axes of the torus $T^2$ with positive orientation, so the index is 1. (Note that the angular coordinate of the principal directions is chosen in the interval $[-1/4,1/4]$, corresponding to principal vectors $(\xi, \eta)$ with $\xi \ge 0$. The jumps in value of $\psi$ from $1/4$ to $-1/4$ are not discontinuities in the line field.) The right panel plots 40 lines of curvature that intersect the blue circle in equally spaced points. This also illustrates that the index is 1 because all of the lines of curvature are transverse to the circle.}
\label{umbi_fig_1}
\end{figure}

For fixed $(x,y)$, the quadratic formula gives non-zero solutions of $R(x,y,\cos(2\pi\psi),\sin(2\pi\psi)) = 0$. Computation of normal curvatures in orthogonal coordinates aligned with the principal directions then determines which of the two principal directions has larger normal curvature. Results from these calculations are displayed in the figures of this section. The figures have pairs of plots which display the following information:
\begin{itemize}
\item
The right plots show numerical integration of forty lines of curvature (black) projected onto the $(x,y)$ plane. The initial points for these computations lie on blue circles of radius $0.1$.                                                                                                                                                                                                                                                                                                                                                     \item 
The left plots show the function $\psi(2\pi\theta)$ with values in the interval $[-1/4,1/4]$, calculated along the circles displayed in the right plots. The orientation of the principal vector $(\xi, \eta)$ is chosen so that $\xi \ge 0$. The colors in right and left plots match. The partitions of the circles were refined so that gaps between adjacent values of $\psi$ are small. 
\end{itemize}
 
Figure~\ref{umbi_fig_1} depicts the principal foliation for $a = (-1,-1,-1,0,0,0,-1,0,-1)$. The left panel plots solutions of $R(0.1\cos(2\pi\theta),0.1\sin(2\pi\theta),\cos(2\pi\psi),\sin(2\pi\psi)) = 0$ on the circle of radius $0.1$ showing that the index is 1. The right hand panel plots 40 lines of curvature intersecting the blue circle at equally spaced points. The eight values of $\theta = \psi$ in the left panel approximate directions in which lines of curvature approach the origin. Note that all the lines of curvature intersect the blue circle transversally, again demonstrating that the index is 1.

\begin{figure}[ht]
    \centering
    \begin{subfigure}{6in}
       \centering
            \includegraphics[height=1.75in]{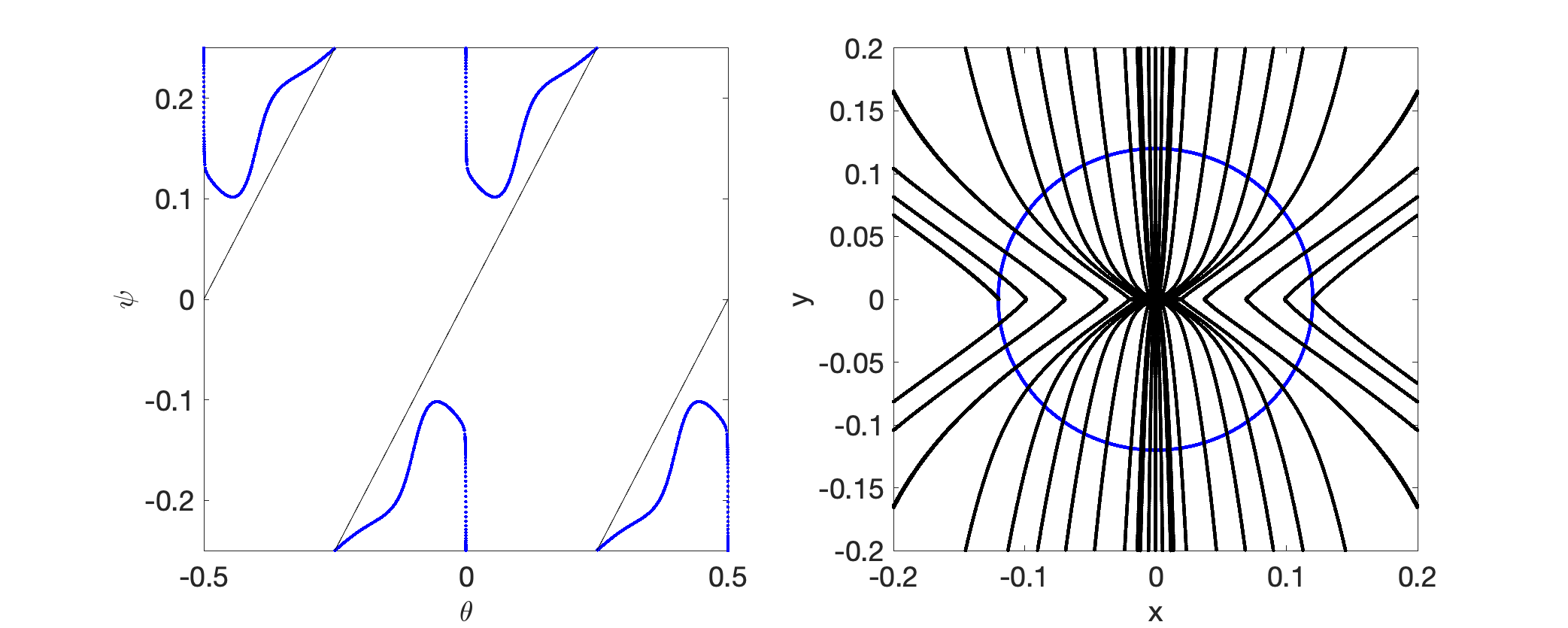}
        \caption{$a = (1/28,0,1,0,0,0,2,0,1)$: index 0}
        \label{umbi_fig_2a}       
    \end{subfigure}
    \begin{subfigure}{6in}
     \centering
            \includegraphics[height=1.75in]{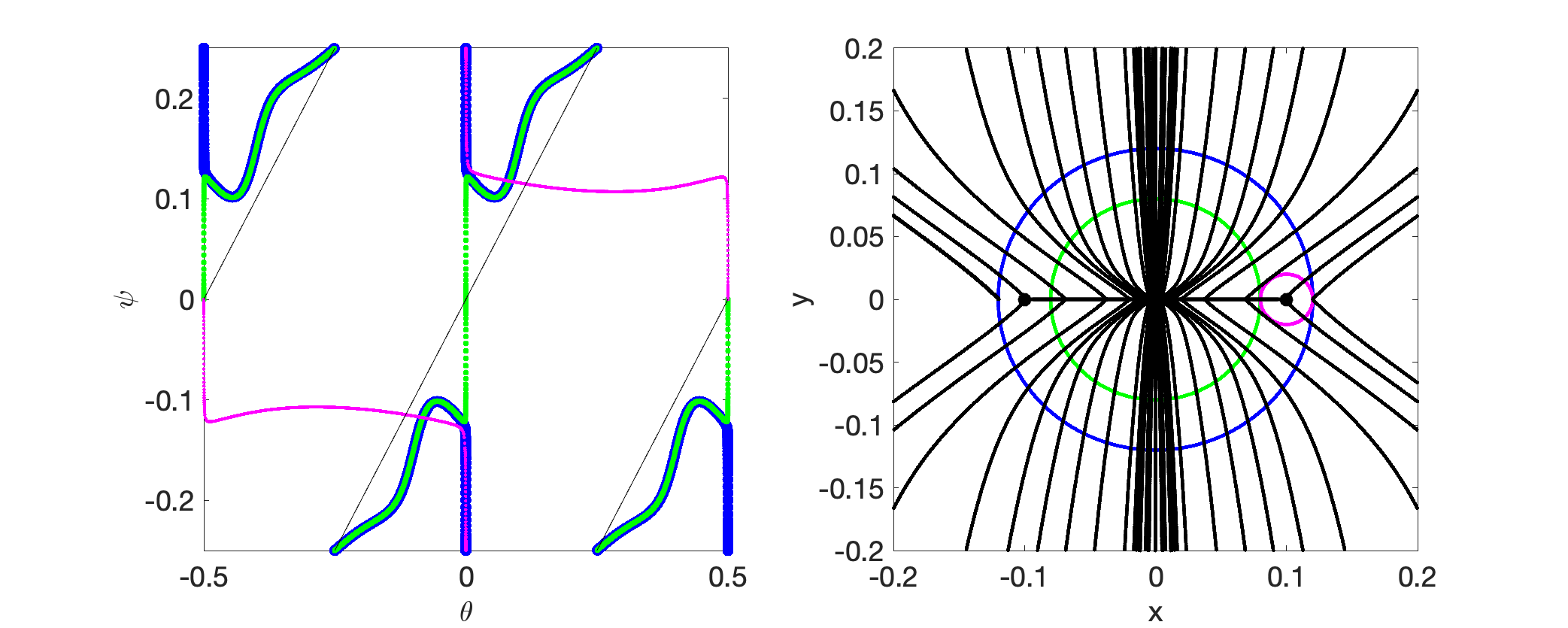}
        \caption{$a = (1/28,0,0.987143,0,0,0,2,0,1)$: index 1}
        \label{umbi_fig_2b}
    \end{subfigure}
    \caption{The index of the umbilic point at the origin jumps as umbilic points merge with the origin. The upper plots have $a_0 = 1/28, \, a_2 = 1$ where a merger occurs.  When  $a_2 < 1$ there  two star umbilics on the $x$ axis. These are located at $x = \pm0.1$ and marked by large black dots in the lower right figure. The left plots are graphs of angular coordinates of principal vectors along the circles displayed in the right plots. The winding number along the blue curves is $0$ in both upper and lower plots, as evident that the graphs are not onto. The winding number of the green curve is $1$, and that is the index of the origin in the lower right plot. The winding number of the  magenta curve is $-1/2$, the index of the star umbilic at $(x,y) = (0.1,0)$}
    \label{umbi_fig_2}
\end{figure}

Figure~\ref{umbi_fig_2} illustrates the change of index that occurs when a pair of umbilic points merge with the origin. The results of the previous section demonstrate that two star umbilic points emerge from the origin as $a_2$ decreases from 1 at parameters $a = (1/28,0,1,0,0,0,2,0,1)$. Figure~\ref{umbi_fig_2a} displays results for $a = (1/28,0,1,0,0,0,2,0,1)$ where the index is $1$. Figure~\ref{umbi_fig_2b} shows results for $a_2 = 0.987143$ where there are two star umbilics near the points $(x,y) = (\pm 0.1,0)$. The left panel of the figure shows the principal directions $\psi \in [-1/4,1/4]$ along two circles centered at the origin of radius $0.12$ (blue) and $0.08$ (green). The winding number around the larger blue circle is 0 while the winding number around the smaller green circle is 1. This is the index of the umbilic at the origin. The small magenta circle in the right hand plot encircles the umbilic point at $0.1$. The winding number of the principal directions around this circle is $-1/2$. Consistent with the index theorem, the winding number of the blue circle is the sum of the indices at the three umbilics in its interior. 

\begin{figure}
    \centering
    \begin{subfigure}[b]{6in}
     \centering
        \includegraphics[height=1.75in]{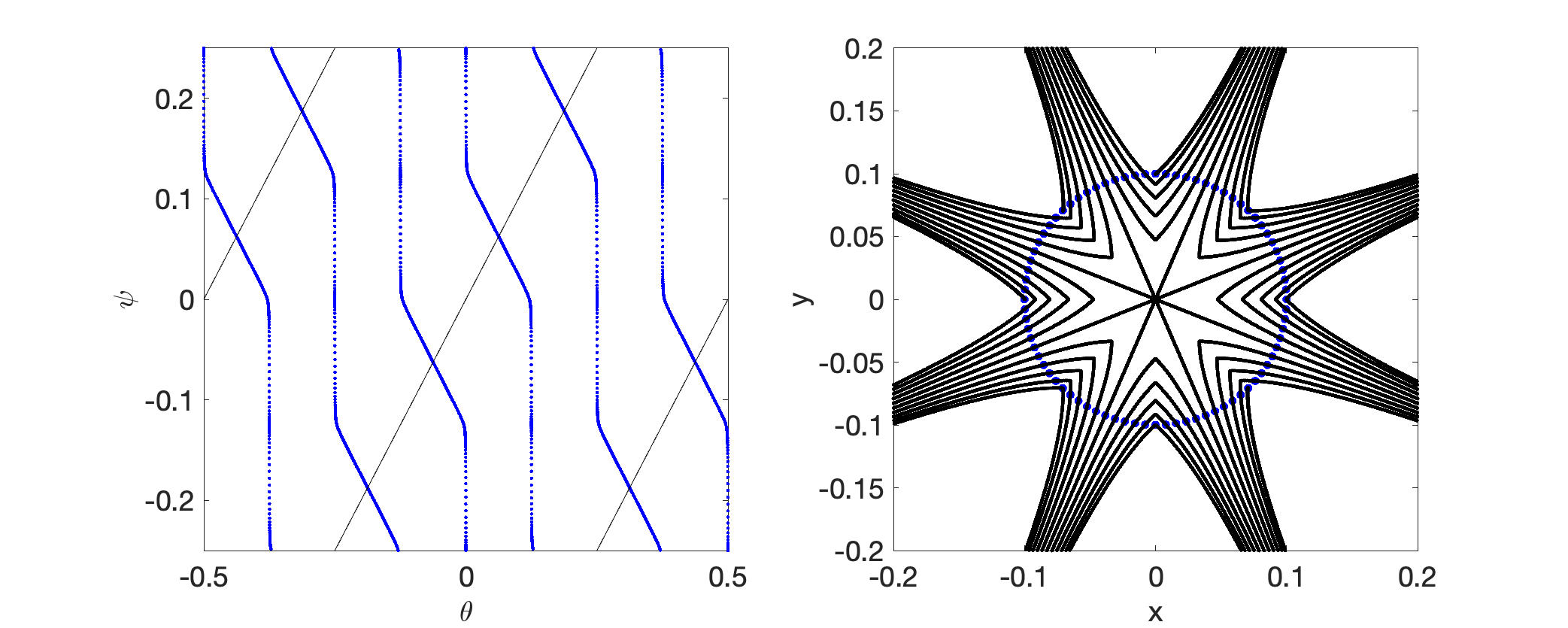}
        \caption{$a = (1/28, 0, 1, 0, -3/2, 0, 1, 0, 1/28)$: index -3}
        \label{umbi_fig_3a}
    \end{subfigure}
    \begin{subfigure}[b]{6in}
     \centering
        \includegraphics[height=1.75in]{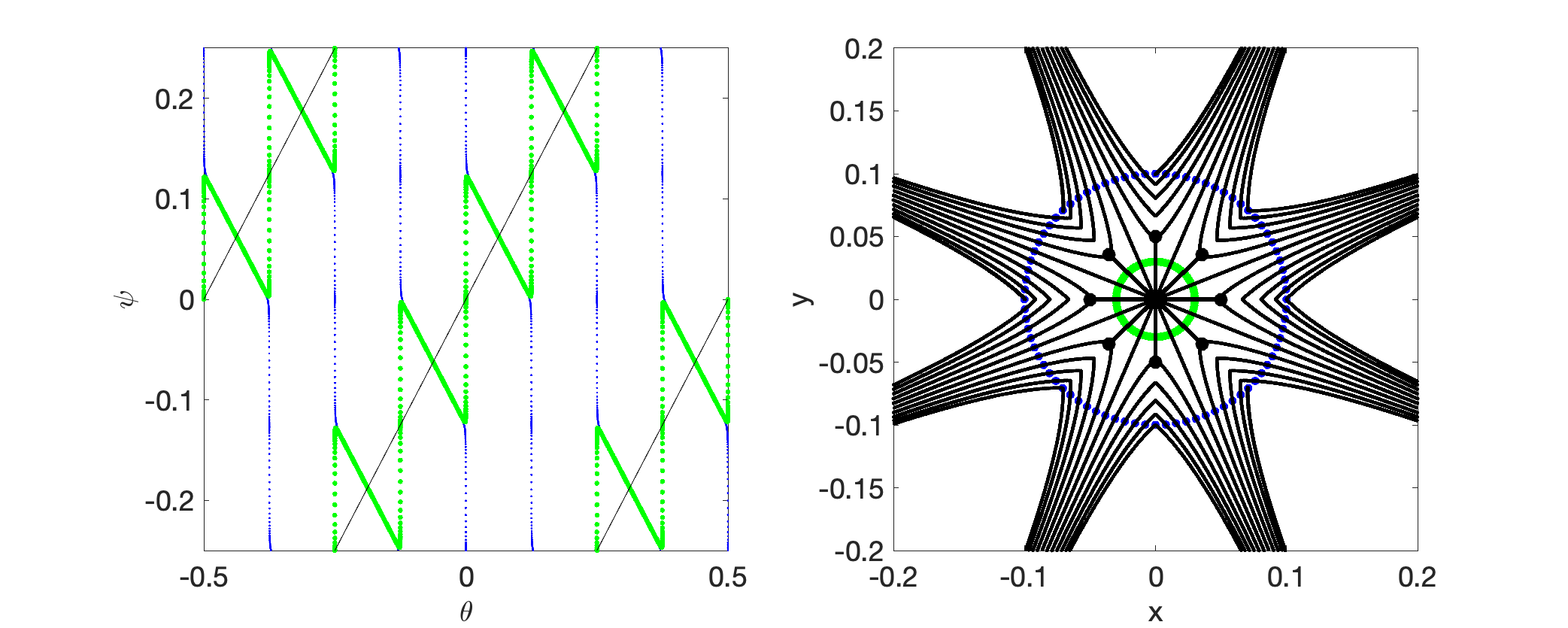}
        \caption{$a = (0.03583,0,1,0,-1.498,0,1,0,0.03583)$: index 1}
        \label{umbi_fig_3b}
    \end{subfigure}
    \caption{Symmetric bifurcation producing jump in index from -3 to +1 at the origin. The angular coordinates in the left hand plots are the same as those in Figure~\ref{umbi_fig_2}.}
    \label{umbi_fig_3}
\end{figure}

The equations for non-zero umbilics on circles of radius r have a limit as $r \to 0$ that are linear in $a$. Since $a$ is a 9 dimensional vector, the maximum number of non-zero umbilics is 8. Figure~\ref{umbi_fig_3} illustrates principal foliations in which surfaces invariant with respect to rotation by $\pi/4$ have 8 umbilics that merge with the origin. 
The merger occurs for parameters $a = (1/28, 0, 1, 0, -3/2, 0, 1, 0, 1/28)$, a polynomial that is invariant with respect to rotation by $\pi/4$. The left panel of Figure~\ref{umbi_fig_3a} shows that the index of the origin is $-3$. The right panel of Figure~\ref{umbi_fig_3b} shows the principal foliation at parameters  $a = (0.03583,0,1,0,-1.498,0,1,0,0.03583)$, chosen to preserve symmetry under rotation by $\pi/4$ and to have 8 umbilics near the circle of radius $0.05$ in the $(x,y)$ plane. These umbilics, depicted by large black dots, are stars with one of their separatrices connecting to the origin along radial lines in the $(x,y)$ plane. Figure~\ref{umbi_fig_3b} plots principal directions along two circles: the larger blue circle has radius $0.1$ and the smaller green circle has radius $0.03$. The winding number $1$ of the principal foliation around the green circle is the index of the origin since the origin is the only umbilic in its interior. The winding number $-3$ of the principal foliation around the blue circle is unchanged from the parameters used for Figure~\ref{umbi_fig_3a}. The principal directions at both sets of parameter values are closely aligned at angles that are not integer multiples of $\pi/4$. At those exceptional angles, the angular coordinates of the blue and green principal directions jump by $1/4$ in opposite directions.

\section*{Centers and Symmetry}

The index of a generic umbilic is $\pm\frac{1}{2}$. Umbilics with index $1$ must be degenerate. Moreover, they are not found in generic one parameter families of surfaces~\cite{MR2105778}. Here we discuss briefly how umbilics with index $1$ appear in generic two parameter families of surfaces and the qualitative changes in principal foliations that occur when these surfaces are perturbed. This issue has hardly been discussed in the literature about principal foliations since umbilics with index $1$ are not found in generic one parameter families of surfaces. Though fragmentary, our results present a conjectural picture of how principal foliations bifurcate in generic two parameter families of surfaces near an umbilic of index 1. 

We start with an analysis of paraboloids $S_{(a,b,c)}$ in Monge form; i.e. the graphs of functions
$$h_{(a,b,c)}(x,y) = ax^2+2bxy+cy^2$$
with $b^2 - ac < 0$. When $a=c$ and $b=0$, $S_{(a,0,a)}$ is a surface of revolution, so one of its principal foliations consists of intersections with planes contain the $z$ axis of symmetry and the second principal foliation consist of circles that are intersections of $S_{(a,0,a)}$ with planes orthogonal to the z-axis. When $b \ne 0$, rotation of the $(x,y)$ plane by angle $\frac{1}{2} \arctan(\frac{2b}{c-a})$ eliminates the middle term $xy$ in the function $h$. This reduces the search for umbilics to the case $b=0$. Calculation of the Rodrigues determinant for the zero level set of $h(x,y,z) = z - ax^2 - cy^2, \, a,c> 0$ finds umbilics at $(\pm ((\frac{c-a}{4a^2c})^{1/2},0)$ when $a < c$ and at  $(0, \pm ((\frac{a-c}{4ac^2})^{1/2})$ when $c < a$. As $a-c$ changes sign from negative to positive, a pair of umbilics on the $x$ axis coalesce at the origin and then reemerge along the $y$ axis. The geometry resembles that of complex square roots of a real variable. In two dimensional families of surfaces  $S_{(a,b,c)}$ with fixed $\frac{a+c}{2}$ and varying $b$ and $\frac{a-c}{2}$, the surface swept out by umbilics is a double cover of this parameter plane which is ramified at the origin. 
 
\begin{figure}
    \centering
    \begin{subfigure}[b]{6in}
        \includegraphics[width=6in,height=1.5in]{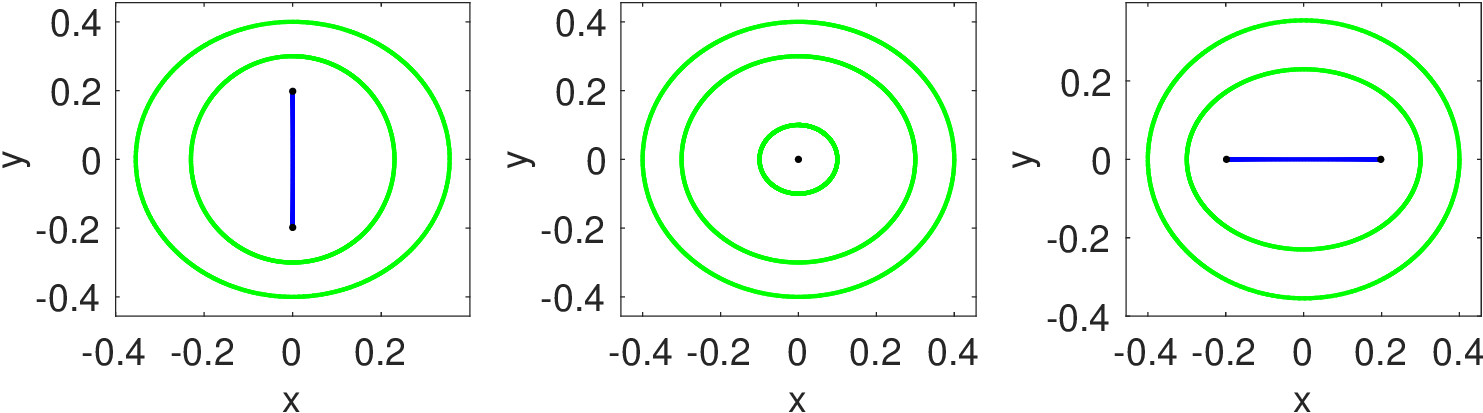}
        \caption{Principal foliations of the graphs of $h(x,y) = x^2/2 + y^2/2 - \lambda(x^2-y^2)$ for (left) $\lambda = -0.01$, (center) $\lambda = 0$ and (right) $\lambda = 0.01$. Apart from umbilics plotted as black dots and connecting lines of curvature plotted blue, all other lines of curvature are closed curves plotted green. As $\lambda$ increases through $0$, the two umbilics on the $y$-axis merge at the origin and remerge along the $x$-axis.}
        \label{umb_fig_6a}
    \end{subfigure}
    \begin{subfigure}[b]{6in}
        \includegraphics[width=6in,height=1.5in]{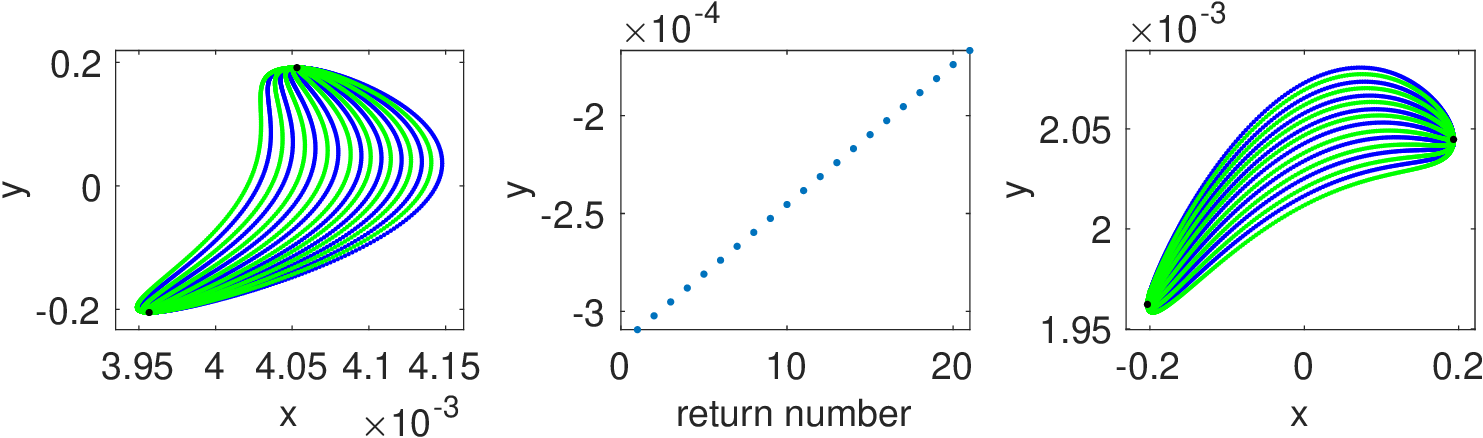}
        \caption{Principal foliation of the graphs of $h(x,y) = x^2/2 + y^2/2 + \lambda(x^2-y^2) -0.002x^3 + 0.002x^2y + 0.001xy^2 -0.001y^3$ for (left) $\lambda = -0.01$ and (right) $\lambda = 0.01$. The left and right panels plot the umbilics together with their separatrix lines of curvature. These separatrices no longer connect umbilics as in Figure~\ref{umb_fig_6a}, instead making tight turns around the opposite umbilic and then spiraling outward. The center panel plots successive values of $y$ for intersections of a single line of curvature with the $x$-axis to the left of the origin at parameter $\lambda = 0$. The line of curvature is not closed: it spirals around the origin slowly.}
                  \label{umb_fig_6b}
    \end{subfigure}
    \caption{Bifurcations of a principal foliation at a degenerate umbilic of index $1$ that splits into two lemon umbilics of index $\frac{1}{2}$.}
    \label{fig:umbilic_split}
\end{figure}

The surfaces $S_{(a,0,c)}$ are symmetric with respect to reflections along the $x$ and $y$ axes since  $h_{(a,0,c)}(x,y) = h_{(a,0,c)}(-x,y)$ and $h_{(a,0,c)}(x,y) = h_{(a,0,c)}(x,-y)$. In the case $0<a<c$, the umbilic points lie in the plane $y=0$ and the intersection of $S_{(a,0,c)}$ with this plane contains a line of curvature connecting the two umbilic points. Furthermore, this connecting line of curvature is surrounded by a family of closed lines of curvature that are symmetric with respect to the reflection $h_{(a,0,c)}(x,y) = h_{(a,0,c)}(x,-y)$. Figure~\ref{umb_fig_6a} displays principal foliations of these paraboloids. Applying rotations as described above, all of the surfaces $S_{(a,b,c)}$, $(a,b,c) \ne (0,0,0)$ have a line of curvature connecting its two umbilics, surrounded by a family of closed lines of curvature. These geometric properties are not generic within the space of smooth, analytic or even polynomial surfaces. 

Generic surfaces lack lines of curvature connecting umbilics and have isolated closed lines of curvature~\cite{MR724448,MR2007065}. Figure~\ref{umb_fig_6b} depicts principal foliations of three surfaces belonging to a one parameter family of surfaces in Monge form obtained by adding the cubic term $-0.002x^3 +0.002x^2y+0.001xy^2-0.001y^3$ to $h_{(1/2 + \lambda,0,1/2-\lambda)}(x,y)$. The parameter $\lambda$ takes the values $-0.01$ in the left panel of the figure and the umbilics lie close to the $y$-axis. (Note the difference in scales of the $x$ and $y$ axes.) The separatrices of the two umbilics are disjoint and plotted as green and blue curves. Each makes a tight turn around the opposite umbilic and spirals outward. The right panel of Figure~\ref{umb_fig_6b} is similar to the left with $\lambda = 0.01$ and the roles of the $x$ and $y$ axes interchanged. The middle panel displays twenty consecutive intersections of a single line of curvature with the ray $x=0;\, y<0$ on the surface with parameter $\lambda = 0$. In contrast to the principal foliation in the center panel of Figure~\ref{umb_fig_6a}, this line of curvature is not closed: it spirals slowly in the radial direction. A more complete bifurcation analysis of surfaces with an umbilic of index $1$ would characterize which surfaces have lines of curvature that connect umbilics and surfaces where the number of closed lines of curvature changes. That analysis is complementary to the main result of this paper related to the Carath\'{e}odory conjecture and not pursued further here.

\section*{Appendix: Prior Work on the Carath\'{e}odory conjecture}
The papers listed below investigate the Carath\'{e}odory conjecture. The papers of Bol, Klotz, Scherbel and Ivanov give successive attempts to complete the arguments initiated by Hamburger. Smyth, Xavier and Guilfoyle use PDE methods and discuss smoothness extensively. Lasarovici and Smyth approach the problem topologically by formulating results that rule out ``elliptic sectors'' of an umbilic. The papers of Ando, Bates, Ghomi and Guilfoyle present ``counterexamples'' with relaxed smoothness or perturbations of the Euclidean metric. This paper extends results of Ando that analyze the index of umbilics at the origin in surfaces that are graphs of homogeneous polynomials.
\begin{itemize}
\item
N. Ando, Umbilics of surfaces that are graphs of homogeneous polynomials, 2000-2002 \cite{MR1789060,MR1955625}
\item
N. Ando, F. Toshifumi and M. Umehara, $C^1$ umbilics of high index, \cite{MR3667760}
\item
 L. Bates, A weak counterexample to the Carath\'{e}odory conjecture, 2001 \cite{MR1845177}
\item
 G. Bol, {\"U}ber Nabelpunkte auf einer Eifl{\"a}che, 1943 \cite{Bol1943}
 \item
 R. Garcia andC. Gutierrez, Ovaloids of {${\bf R}^3$} and their umbilics: a differential equation approach, 1998, \cite{MR1801351}
 \item
 M. Ghomi and R. Howard, Normal curvatures of asymptotically constant graphs and Carath\'{e}odory's conjecture, 2012, \cite{MR2957223}
 \item
 B. Guilfoyle, On isolated umbilic points, 2020 \cite{MR4201794} 
\item
 B, Guilfoyle and W. Klingenberg, Isolated umbilical points on surfaces in {$\mathbb R^3$}, 2006 \cite{MR3602286}
 \item
 C. Gutierrez, F. Mercuri  and F. S\'{a}nchez-Bringas, On a conjecture of Carath\'{e}odory: analyticity versus smoothness, 1996 \cite{MR1412952}
 \item
 C. Gutierrez and F. S\'{a}nchez-Bringas, On a Carath\'{e}odory's conjecture on umbilics: representing ovaloids, 1997 \cite{MR1492978}
 \item
 C. Gutierrez and J. Sotomayor, Lines of curvature, umbilic points and Carath\'{e}odory conjecture, 1998 \cite{MR1633013}
\item 
 H. Hamburger, Beweis einer Carath\'{e}odoryschen {V}ermutung I,II,III, 1940-41 \cite{MR1052,MR6480,MR6481}
 \item
 V. Ivanov, An analytic conjecture of Carath\'{e}odory, 2002 \cite{MR1902826}
\item
 T. Klotz, On G. Bol's proof of Carath{\'e}odory's conjecture, 1959 \cite{MR120602}
 \item
 L. Lazarovici, Elliptic sectors in surface theory and the Carath\'{e}odory-{L}oewner conjectures, 2000 \cite{MR1863731}
 \item
 B. Smyth, The nature of elliptic sectors in the principal foliations of surface theory, 2005 \cite{MR2185156}
\item
 B. Smyth and F. Xavier, A sharp geometric estimate for the index of an umbilic on a smooth surface, 1992 \cite{MR1148679}
 \item
 B. Smyth and F. Xavier, Real solvability of the equation {$\partial^2_{\overline z}\omega=\rho g$} and the topology of isolated umbilics, 1998 
 \item
 H. Scherbel, A new proof of {H}amburger's {I}ndex {T}heorem on umbilical points, 1993 \cite{MR2714735}
  \item
 J. Sotomayor and L. Mello, A note on some developments on Carath{\'e}odory conjecture on umbilic points, 1999 \cite{Sotomayor1999}
\item
 C. Titus, A proof of a conjecture of Loewner and of the conjecture of Carath{\'e}odory on umbilic points, 1973 \cite{Titus1973}
  \item
 F. Xavier, An index formula for {L}oewner vector fields, 2007 \cite{MR2350132}
\end{itemize}
In addition to these publications, Guilfoyle posted a manuscript on the ArXiv in 2008, that presents a different strategy for proving the Carath\'{e}odory conjecture. 

Matlab codes for producing the figures in this paper can be found in the Github repository umbilic\_index\_2\_2025.


\end{document}